\documentclass{elsarticle-modified}

\usepackage{amssymb,amsmath,amsfonts}
\usepackage{xcolor}

\journal{Appl. Math. Lett.}

\newtheorem{theorem}{Theorem}

\newtheorem{proposition}{Proposition}
\newtheorem{remark}{Remark}

\newcommand\nE{\mathcal{E}}
\newcommand\nS{\mathcal{S}}

\newcommand{\ee}{\mathrm{e}}

\newcommand{\Id}{\mathrm{Id}}

\begin{document}

\begin{frontmatter}

\title{Non-existence of generalized splitting methods with positive coefficients of order higher than four}

\author{Winfried Auzinger}
\address{Technische Universit{\"a}t Wien,
Institut f{\"u}r Analysis und Scientific Computing,
Wiedner Hauptstrasse 8--10/E101, A-1040 Wien, Austria}
\ead{w.auzinger@tuwien.ac.at}
\ead[url]{www.asc.tuwien.ac.at/~winfried/}

\author{Harald Hofst\"atter\corref{ourcorrespondingauthor}}
\cortext[ourcorrespondingauthor]{Corresponding author}
\address{Universit{\"a}t Wien, Institut f{\"u}r Mathematik, Oskar-Morgenstern-Platz 1, A-1090 Wien, Austria}
\ead{hofi@harald-hofstaetter.at}
\ead[url]{www.harald-hofstaetter.at}

\author{Othmar Koch}
\address{Universit{\"a}t Wien, Institut f{\"u}r Mathematik, Oskar-Morgenstern-Platz 1, A-1090 Wien, Austria}
\ead{othmar@othmar-koch.org}
\ead[url]{www.othmar-koch.org}

\begin{abstract}
We prove that generalized exponential splitting methods making explicit use of commutators
of the vector fields are limited to order four when only real coefficients are admitted.
This generalizes the restriction to order two for classical splitting methods with only positive coefficients.
\end{abstract}


\end{frontmatter}



\section{Introduction}
\subsection{Splitting methods}
We consider evolution equations on $\mathbb{R}^d$ or $\mathbb{C}^d$ where the right-hand side is split into two components,
\begin{equation}\label{eq:splitevo}
\partial_t y(t) = A(y(t))+B(y(t)),\quad t\geq t_0,\quad y(t_0)=y_0.
\end{equation}
In this introduction we only consider the linear case where $A,B$ are linear
operators represented by real or complex matrices.
For the numerical solution of~\eqref{eq:splitevo} we consider $s$-stage  splitting methods, where
one step $(t_n, y_n)\mapsto (t_{n+1},y_{n+1})$ with step-size $\tau$  is given by
\begin{equation}\label{eq:splitting_method}
y_{n+1}=\nS(\tau)y_n=\ee^{b_s\tau B}\ee^{a_s\tau A}\cdots\ee^{b_1\tau B}\ee^{a_1\tau A}y_n.
\end{equation}
A splitting method has convergence order $p$ if it holds
\begin{equation*}
\nS(\tau)y_n=\ee^{\tau(A+B)}y_n+O(\tau^{p+1}).
\end{equation*}
\subsection{Positive coefficients}
For certain applications only splitting schemes with non-negative coefficients are
suitable. In particular this is the case if $A$ is a discretized
sectorial operator associated with a parabolic equation, because in this
case the flow of $A$ is non-reversible and does not tolerate negative
time increments in the numerical approximation, whence $a_j$ is required to be non-negative.
A splitting method of order $p=2$ with all coefficients positive is given by Strang
splitting
\begin{equation*}
\nS(\tau) = \ee^{\frac{1}{2}\tau A}\ee^{\tau B}\ee^{\frac{1}{2}\tau A}.
\end{equation*}
It is known that $p=2$ is the maximum order of a splitting method
with all coefficients positive, a fact which is established by the following theorem.

\begin{theorem}
If $\nS$ is a splitting method~\eqref{eq:splitting_method} of order $p\geq 3$ with real coefficients, then at least one of the coefficients $a_j$
is strictly negative, and also at least one of the
coefficients $b_j$ is strictly negative.
\end{theorem}
This theorem was first proved in~\cite{golkap96}, see also~\cite{BLANES200523}.
A weaker version stating that at least one of all coefficients $a_j$, $b_j$ combined is strictly negative,
was proved earlier in~\cite{sheng86}.

\subsection{Generalized splitting methods}
In many applications the commutator $[B,[B,A]]$ and its exponential are
readily computable, see~\cite{omelyanetal02}. This suggests to consider
generalized splitting methods of the form
\begin{equation}\label{eq:generalized_splitting1}
y_{n+1}=\nS(\tau)y_n=\ee^{c_s\tau^3 [B,[B,A]]}\ee^{b_s\tau B}\ee^{a_s\tau A}\cdots\ee^{c_1\tau^3 [B,[B,A]]}\ee^{b_1\tau B}\ee^{a_1\tau A}y_n,
\end{equation}
or
\begin{equation}\label{eq:generalized_splitting2}
y_{n+1}=\nS(\tau)y_n=\ee^{b_s\tau B+c_s\tau^3 [B,[B,A]]}\ee^{a_s\tau A}\cdots\ee^{b_1\tau B+c_1\tau^3 [B,[B,A]]}\ee^{a_1\tau A}y_n,
\end{equation}
which possibly allow orders higher than 2 while involving only positive coefficients.
Indeed, the scheme
\begin{equation*}
\nS(\tau)=\ee^{\frac{1}{6}\tau B}\ee^{\frac{1}{2}\tau A}
\ee^{\frac{2}{3}\tau B-\frac{1}{72}\tau^3[B,[B,A]]}
\ee^{\frac{1}{2}\tau A}\ee^{\frac{1}{6}\tau B}
\end{equation*}
proposed in~\cite{suzuki95,chin97} has order $p=4$ and positive
coefficients $a_j$ and $b_j$.
However, as established by the following theorem, $p=4$ is the maximum order of such a generalized splitting method
with all coefficients $a_j$ positive. This holds even under the
additional assumption $[B,[B,[B,A]]]=0$, which in many applications is
satisfied, see~\cite{omelyanetal02}.
\begin{theorem}\label{thm:Thm_pos_gen_split}$~$

\medskip\noindent
(i) If $\nS$ is a generalized splitting method
of the form~\eqref{eq:generalized_splitting1}
or~\eqref{eq:generalized_splitting2}
of order $p\geq 5$ with real coefficients, then at least one of the coefficients $a_j$
is strictly negative.

\medskip\noindent
(ii) If $\nS$ is a generalized splitting method~\eqref{eq:generalized_splitting2} with
real coefficients which is of order $p\geq 5$
if applied to an equation~\eqref{eq:splitevo} where the operators $A,B$ satisfy
$[B,[B,[B,A]]]=0$, then at least one of the coefficients $a_j$ is strictly negative.
\end{theorem}
It is clear\footnote{By logical transposition: If an object (here a generalized
splitting method with all coefficients $a_j$ nonnegative) does not exist under
some restrictive assumptions, then it cannot exist under more general assumptions.
\label{fn:contraposition}}
that part (i) follows  immediately from part (ii), which immediately
follows from the following theorem, which may be interesting in itself.
\begin{theorem}\label{thm:Thm_pos_split_special}
If $\nS$ is a  splitting method~\eqref{eq:splitting_method} with
real coefficients which is of order $p\geq 5$
if applied to an equation~\eqref{eq:splitevo} where the operators $A,B$ satisfy
$[B,[B,A]]=0$, then at least one of the coefficients $a_j$ is strictly negative.
\end{theorem}
A proof of Theorem~\ref{thm:Thm_pos_gen_split} was proposed in~\cite{chin2005}.
In Section~\ref{Sect:proof_of_thm3}
we will give a new independent proof by showing that
Theorem~\ref{thm:Thm_pos_split_special}
(and thus also Theorem~\ref{thm:Thm_pos_gen_split}) is an easy consequence of
a recent result proved by the authors in~\cite{positivitypaper2018}.

\section{Proof of Theorem~\ref{thm:Thm_pos_split_special}}
\label{Sect:proof_of_thm3}
The essential step leading to the main result of~\cite{positivitypaper2018} is comprised by the following
proposition.\footnote{Note that we have changed some denotations: $H(t)$, $H_0$, $H_1$, $s$, $a_j$, $c_j$ correspond
respectively to the denotations $A(t)$, $A_0$, $A_1$, $J$, $b_j$, $y_j$ of
\cite{positivitypaper2018}.}
\begin{proposition}\label{prop:from_pos_paper}
Let $u(t)$ be the exact solution of
\begin{equation}\label{eq:nonauto_special}
\partial_t u(t) = H(t)u(t) = (H_0+tH_1)u(t),\quad u(0)=u_0
\end{equation}
with $H_0, H_1\in\mathbb{C}^{d\times d}$
and
\begin{equation}\label{eq:comm_free_step}
v(\tau)=\ee^{a_s\tau H_0+c_s\tau^2H_1}
\cdots\ee^{a_1\tau H_0+c_1\tau^2H_1}u_0
\end{equation}
with given coefficients $a_j,c_j\in\mathbb{R}$.
If
\begin{equation*}
v(\tau)-u(\tau)=O(\tau^6),
\end{equation*}
then at least one of the coefficients $a_j$ is
strictly negative.\footnote{{See Remark~\ref{rem:strictly_negative} below.}}
\end{proposition}
Here~\eqref{eq:comm_free_step} can be interpreted as one step with step-size $\tau$
of a commutator-free exponential integrator applied
to the special non-autonomous equation~\eqref{eq:nonauto_special}.
To show that Theorem~\ref{thm:Thm_pos_split_special} follows
from Proposition~\ref{prop:from_pos_paper} we
first  use the standard reformulation
\begin{equation}\label{eq:auton_reformulation}
y(t)=\left({s(t)\atop u(t)}\right),\quad
\partial_t y(t)=\underbrace{\left({0\atop H(s(t))u(t)}\right)}_{A(y(t))}
+\underbrace{\left({1\atop 0}\right)}_{B(y(t))}, \quad
y(0)=y_0=\left({0 \atop u_0}\right)
\end{equation}
of the non-autonomous problem~\eqref{eq:nonauto_special}
as an autonomous problem by adding the component $ s(t)=t $. Here the operators
$A,B\colon \mathbb{R}^{d+1}\to\mathbb{R}^{d+1}$ are nonlinear, therefore
a direct application of the splitting method~\eqref{eq:splitting_method} is not
possible. However, by associating  the flows $\nE_A(t, y_0)$, $\nE_B(t, y_0)$
of the subproblems $\partial_t y(t)=A(y(t))$, $\partial_t y(t)=B(y(t))$
with exponentials of
Lie derivatives\footnote{We adopt the notation
                         from~\cite[Chapter~III]{haireretal02b}.}
$\ee^{t D_A}$, $\ee^{t D_B}$, which act on a smooth
map $F:\mathbb{R}^{d+1}\to\mathbb{R}^{d+1}$ as
\begin{equation*}
(\ee^{t D_A}F)(y) = F(\nE_A(t, y)),
\qquad
(\ee^{t D_B}F)(y) = F(\nE_B(t, y)),
\end{equation*}
and thus
\begin{equation*}
\nE_A(t, y_0) = (\ee^{t D_A}\Id)(y_0),\quad
\nE_B(t, y_0) = (\ee^{t D_B}\Id)(y_0),
\end{equation*}
each splitting method~\eqref{eq:splitting_method} of order $p$
for linear problems~\eqref{eq:splitevo} can be promoted to a splitting method
\begin{align}
y_{n+1}=
  \nS(\tau, y_n)&=\nE_B(b_s\tau,\cdot)\circ\nE_A(a_s\tau,\cdot)
  \circ\ldots\circ\nE_B(b_1\tau,\cdot)\circ\nE_A(a_1\tau,y_n)  \nonumber\\
  &=\big(\ee^{a_1\tau D_A}\ee^{b_1\tau D_B}\cdots\ee^{a_s\tau D_A}\ee^{b_s\tau D_B}
  \Id\big)(y_n) \label{eq:splitting_method2}
\end{align}
of the same order for nonlinear problems, see~\cite[Section~III.5.1]{haireretal02b}.

\begin{remark}\normalfont\label{rem:abstract_lie}
 The convergence order of a (generalized) splitting method is
determined by order conditions, which are polynomial equations in the coefficients
of the method. Usually these conditions are derived in a purely formal way in the abstract
algebra of formal power series in the non-commuting variables $\mathtt{A}$, $\mathtt{B}$
and its embedded Lie algebra with Lie bracket defined by $[X,Y]=XY-YX$,
see~\cite{auzingeretal13c,MuntheKaas957}.
By associating  $\mathtt{A}$, $\mathtt{B}$ with the matrices  $A$, $B$ in the linear case, and with the Lie derivatives
$D_A$, $D_B$ in the nonlinear case, it follows
that~\eqref{eq:splitting_method} and~\eqref{eq:splitting_method2} indeed have
the same order~\cite{haireretal02b}.
\end{remark}

For the special problem~\eqref{eq:auton_reformulation}
the Lie derivatives are given by
\begin{equation*}
D_A=\sum_{i=1}^{d+1}A_i(y)\frac{\partial}{\partial y_i}
=\sum_{i=1}^d\,{[H(s)u]}_i\frac{\partial}{\partial u_i}
=\sum_{i=1}^d\,{[(H_0+sH_1)u]}_i\frac{\partial}{\partial u_i}
\end{equation*}
and
\begin{equation*}
D_B=\sum_{i=1}^{d+1}B_i(y)\frac{\partial}{\partial y_i}
=\frac{\partial}{\partial s}.
\end{equation*}
A straightforward calculation leads to
\begin{equation*}
[D_B,D_A] = \sum_{i=1}^d\,{[H'(s)u]}_i\frac{\partial}{\partial u_i}
=\sum_{i=1}^d\,{[H_1u]}_i\frac{\partial}{\partial u_i}
\end{equation*}
and
\begin{equation*}
[D_B,[D_B,D_A]] = \sum_{i=1}^d\,{[H''(s)u]}_i\frac{\partial}{\partial u_i} = 0,
\end{equation*}
which shows that the condition $[B,[B,A]]=0$  of
Theorem~\ref{thm:Thm_pos_split_special} promoted to the nonlinear case is satisfied.
Repeated application of the formal identity
\begin{equation*}
\ee^X\ee^Y=\ee^{Y+[X,Y]]}\ee^X\ \ \mbox{for}\ \ [X,[X,Y]]=0
\end{equation*}
to~\eqref{eq:splitting_method2} yields
\begin{align}\label{eq:splitting_transformed}
&\nS(\tau,y_0) = {} \\
& {} = \big(\ee^{a_1\tau D_A}\ee^{a_2\tau D_A+c_2\tau^2[D_B,D_A]}\cdots
            \ee^{a_s\tau D_A+c_s\tau^2[D_B,D_A]}\ee^{(b_1+\ldots+b_s)\tau D_B}\Id\big)(y_0)
   \notag
\end{align}
with well-defined coefficients $c_2,\dots,c_s\in\mathbb{R}$.
Here a single exponential acts as
\begin{align*}
\big(\ee^{a_j\tau D_A+c_j\tau^2[D_B,D_A]}F\big)\left({s\atop u}\right)
&=F\left({s \atop \ee^{a_j\tau H(s)+c_j\tau^2 H'(s)}u}\right) \\
&=F\left({s \atop \ee^{a_j\tau H_0+a_j\tau s H_1+c_j\tau^2 H_1}u}\right),
\end{align*}
and thus, substituting $s=0$,
\begin{equation*}
\big(\ee^{a_j\tau D_A+c_j\tau^2[D_B,D_A]}F\big)\left({0\atop u}\right)
= F\left({0 \atop \ee^{a_j\tau H_0+c_j\tau^2 H_1}u}\right).
\end{equation*}
It follows that for $y_0=\left({0\atop u_0}\right)$ the lower components
of~\eqref{eq:splitting_transformed} can be written as
\begin{equation*}
\ee^{a_s\tau H_0+c_s\tau^2 H_1}\cdots\ee^{a_2\tau H_0+c_2\tau^2 H_1}\ee^{a_1\tau H_0}u_0,
\end{equation*}
which is of the form~\eqref{eq:comm_free_step}.
From Proposition~\ref{prop:from_pos_paper} it follows
that if the splitting method~\eqref{eq:splitting_method2} has order $p\geq 5$
if applied to the special problem~\eqref{eq:auton_reformulation},
or, a fortiori\footnote{See footnote~\ref{fn:contraposition}.},
if applied to nonlinear problems with $[B,[B,A]]=0$ in general,
then at least one of the coefficients $a_j$ is
strictly negative. We have thus proved the nonlinear version
of Theorem~\ref{thm:Thm_pos_split_special}.
For similar formal reasons as in Remark~\ref{rem:abstract_lie},
the linear version of Theorem~\ref{thm:Thm_pos_split_special} follows as well.
\begin{remark}\label{rem:strictly_negative}
Strictly speaking, only a version of
Proposition~\ref{prop:from_pos_paper} with the weaker conclusion that
at least one of the coefficients $a_j$ is  {\em non-positive} has been
proved in~\cite{positivitypaper2018}. Since we may assume form the outset
that $a_j\neq 0$ for $j=2,\dots,s$ in~\eqref{eq:splitting_method},
it is clear that this weaker version already suffices for the proof of
Theorem~\ref{thm:Thm_pos_split_special}. Conversely,
Proposition~\ref{prop:from_pos_paper}  follows from
Theorem~\ref{thm:Thm_pos_split_special}, as can be shown
by a similar reasoning as before.
Thus, the version of Proposition~\ref{prop:from_pos_paper}  given here
follows from the weaker version proved in~\cite{positivitypaper2018}.
\end{remark}
\section*{Acknowledgements}

This work was supported in part by the Vienna Science and Technology Fund (WWTF) [grant number MA14-002].


\end{document}